%% file: manuscript.tex
\documentclass[final,3p,times,fleqn]{elsarticle}

\usepackage{amsmath,amssymb,amsfonts,amscd,mathtools}

\usepackage{algorithm}
\usepackage{algpseudocode}
\usepackage{graphicx}
\usepackage{textcomp}
\usepackage{xcolor}
\usepackage{lineno}
\modulolinenumbers[5]
\usepackage[colorlinks=true, allcolors=blue]{hyperref}
\usepackage{cleveref}

\usepackage{subcaption}
\usepackage[textsize=small,linecolor=red,backgroundcolor=red!25,bordercolor=red]
  {todonotes}

\usepackage{booktabs,cellspace}
\usepackage{siunitx}

\usepackage[toc,page]{appendix}

\renewcommand{\vec}[1]{\boldsymbol#1} 
\usepackage{stmaryrd}

\usepackage{tikz}
\usepackage{pgfplots}

\definecolor{mycolor1}{rgb}{0.00000,0.44700,0.74100}%
\definecolor{mycolor2}{rgb}{0.85000,0.32500,0.09800}%
\definecolor{mycolor3}{rgb}{0.92900,0.69400,0.12500}%
\definecolor{mycolor4}{rgb}{0.49400,0.18400,0.55600}%
\definecolor{mycolor5}{rgb}{0.46600,0.67400,0.18800}%
\definecolor{mycolor6}{rgb}{0.30100,0.74500,0.93300}%
\definecolor{mycolor7}{rgb}{0.63500,0.07800,0.18400}%

\newcommand{\Red}[1]{\textcolor{black}{#1}}
\newcommand{\Blue}[1]{\textcolor{black}{#1}}
\pgfplotsset{
        colormap={blueToRed}{
            rgb255=(0,0,255)
            rgb255=(255,255,255)
            rgb255=(255,0,0)
        },
    }
    
\pgfplotsset{
        colormap={coolToWarm}{
    rgb255=( 59, 76,192)
    rgb255=( 97,128,233)
    rgb255=( 97,128,233)
    rgb255=(119,154,247)
    rgb255=(141,177,254)
    rgb255=(166,197,254)
    rgb255=(191,212,245)
    rgb255=(221,221,221)
    rgb255=(245,196,172)
    rgb255=(245,154,123)
    rgb255=(230,113, 88)
    rgb255=(214, 82, 67)
    rgb255=(201, 58, 54)
    rgb255=(190, 36, 45)
    rgb255=(180,  4, 38)
        },
    }

\usepackage{listings}
\begin{document}

\begin{frontmatter}
\journal{Computer Methods in Applied Mechanics and Engineering}

\title{
Multigrid reduction preconditioning framework for coupled processes in porous and fractured media}

\author[CASC]  {Quan~M.~Bui\corref{cor1}} \ead{bui9@llnl.gov}
\author[TOTAL] {Fran{\c c}ois~P.~Hamon}   \ead{francois.hamon@total.com}
\author[AEED]  {Nicola~Castelletto}       \ead{castelletto1@llnl.gov}
\author[CASC]  {Daniel~Osei-Kuffuor}      \ead{oseikuffuor1@llnl.gov}
\author[AEED]  {Randolph~R.~Settgast}     \ead{settgast1@llnl.gov}
\author[AEED]  {Joshua~A.~White}          \ead{white230@llnl.gov}

\cortext[cor1]{Corresponding author.}
\address[CASC]{Center for Applied Scientific Computing, Lawrence Livermore National Laboratory, United States} 
\address[TOTAL]{Total E\&P Research and Technology, Houston, TX, United States}
\address[AEED]{Atmospheric, Earth, and Energy Division, Lawrence Livermore National Laboratory, United States}

\input{Abstract}

\end{frontmatter}

\allowdisplaybreaks


\input{Section1_Introduction}

\input{Section3_MGR}
\input{Section4_Linear_Systems}

\input{Section5_Results}

\input{Section6_Conclusions}

\section*{Acknowledgements}
\label{sec::acknow}
Funding for the solver development in this work was provided by LLNL through Laboratory Directed Research and Development (LDRD) Project 18-ERD-027. Simulator development support and the data used for the field example was provided by Total S.A. through the FC-MAELSTROM and TC02223 Projects. Portions of this work were performed under the auspices of the U.S. Department of Energy by Lawrence Livermore National Laboratory under Contract DE-AC52-07-NA27344. FPH completed this work during a visiting scientist appointment at LLNL.

\appendix
\input{Appendix}
\biboptions{sort&compress}
\bibliographystyle{elsarticle-num}
\bibliography{reference}

\end{document}

%% file: Abstract.tex
\begin{abstract}
Many subsurface engineering applications involve tight-coupling between fluid flow, solid deformation, fracturing, and similar processes.  
To better understand the complex interplay of different governing equations, and therefore design efficient and safe operations, numerical simulations are widely used.  
Given the relatively long time-scales of interest, fully-implicit time-stepping schemes are often necessary to avoid time-step stability restrictions.
A major computational bottleneck for these methods, however, is the linear solver.  These systems are extremely large and ill-conditioned.
Because of the wide range of processes and couplings that may be involved---e.g. formation and propagation of fractures, deformation of the solid porous medium, viscous flow of one or more fluids in the pores and fractures, complicated well sources and sinks, etc.---it is difficult to develop general-purpose but scalable linear solver frameworks.
This challenge is further aggravated by the range of different discretization schemes that may be adopted, which have a direct impact on the linear system structure.
To address this obstacle, we describe a flexible \Red{strategy} based on multigrid reduction (MGR) that can produce purely algebraic preconditioners for a wide spectrum of relevant physics and discretizations.
%
\Red{We demonstrate that MGR, guided by physics and theory in block preconditioning, can tackle several distinct and challenging problems, notably: a hybrid discretization of single-phase flow, compositional multiphase flow with complex wells, and hydraulic fracturing simulations.}
Extension to other systems can be handled quite naturally.
We demonstrate the efficiency and scalability of the resulting solvers through numerical examples of difficult, field-scale problems.
\end{abstract}

\begin{keyword}
preconditioning, algebraic multigrid, hydraulic fracturing, compositional flow, mimetic finite difference method
\end{keyword}

%% file: Section1_Introduction.tex
\section{Introduction}
With recent focus on creating a reliable, but sustainable, energy economy and mitigating the impacts of climate change, subsurface engineering is playing an increasingly important role.  
Over the past two decades, significant investment has been made in carbon sequestration, geothermal energy, unconventional oil and gas production, geologic waste disposal, and aquifer management.
These applications involve a number of engineering controls that can be optimized to maximize efficiency while minimizing costs and risks.
Numerical simulations are frequently used to provide physical insight into how various design choices may impact operations.
These simulations range from purely empirical models to physics-driven simulators, with a spectrum of complexity and computational cost.

In this work, we focus specifically on high-fidelity simulations of flow, transport, and geomechanics.  
Each process is modeled by a time-dependent, nonlinear partial differential equation (PDE) enforcing conservation of mass or momentum balance.
Together, they form a coupled system whose dynamics are governed by the interplay of each individual process.
The preferred time-integration strategy for these systems is fully-implicit to accurately capture coupling  while resolving relatively long time-scales of interest.
The downside of this approach is that it leads to very large and ill-conditioned block linear systems that must be solved at each Newton iteration to advance the nonlinear solver in time.  
The use of black-box preconditioning strategies for these systems often shows poor parallel scalability or simply fails to converge.
Thus, our goal is to develop efficient preconditioners targeting quite general block linear systems based on a multigrid reduction framework.
A flexible, algebraic perspective allows us to tackle a wide range of physics and discretizations encountered in subsurface engineering practice within a common framework, recycling optimized code components in an efficient way.
To date, this lack of general-purpose frameworks has strongly limited the size and complexity of models that can be run or the ease of porting scalable methods to different simulators.

In particular, we tackle $\mathcal{B}\times \mathcal{B}$ block linear systems in the generic form
\begin{linenomath}
\begin{align}
A x &= b,
&
A &=
\begin{pmatrix}
  A_{11} & A_{12} & \cdots & A_{1 \mathcal{B}} \\
  A_{21} & A_{22} & \cdots & A_{2 \mathcal{B}} \\
  \vdots & \vdots & \ddots & \vdots \\
  A_{\mathcal{B} 1} & A_{\mathcal{B} 2} & \cdots & A_{\mathcal{B} \mathcal{B}}
\end{pmatrix},
&
x &=
\begin{pmatrix}
  x_1 \\
  x_2 \\
  \vdots \\
  x_{\mathcal{B}}
\end{pmatrix},
&
b &= 
\begin{pmatrix}
  b_1 \\
  b_2 \\
  \vdots \\
  b_{\mathcal{B}}
\end{pmatrix}. \label{eq:linear_system}
\end{align}
\end{linenomath}
Here, $A_{ij}$ with $i, j = 1, 2, \dots, \mathcal{B}$ are block matrices of size $N_i \times N_j$ induced by a physics-based partitioning of the the unknowns collected in $x$ by field type, with $N_i$ (respectively $N_j$) the numbers of degrees of freedom associated with the physical quantities $x_i$ (respectively $x_j$).
%
Much of the work in developing solver strategies for these systems has focused on block preconditioners with approximate Schur complements, usually using physics-based strategies. 
\Blue{Some notable physics-based preconditioners include fixed-stress splitting \cite{White16} for coupled flow and mechanics, and \textit{Constrained Pressure Residual} (CPR) \cite{Wallis83, Wallis85, Cao05} for compositional multiphase flow, including also thermal effects \cite{Roy_JCP_19,Roy_SISC_20,Cremon_JCP_20}.}
This approach uses knowledge of the specific physical processes involved to break the tightly coupled systems into smaller sub-problems whose properties are well-studied and more easily handled. 
%
%
As such, block preconditioners are often highly specialized and domain-specific. 
There have been efforts to facilitate development of more general-purpose block preconditioners by combining existing single-physics solvers. 
For example, PETSc's `PCfieldsplit' functionality allows one to compose sophisticated block preconditioners by splitting the solution vector into physical fields \cite{Baley12}. 
However, for systems with a large number of blocks, this task still requires significant effort as all the components of such preconditioners, i.e. the factorization of the global block system, approximations of the Schur complements, solver choices for each block, etc., needs to be explicitly defined by the user. 
Similarly, the Teko package of the Trilinos project \cite{Heroux05} provides a tool for automated construction of block preconditioners, though mostly tailored to the Navier-Stokes equation.

Here, we present an algebraic framework to target general block problems using \textit{multigrid reduction} (MGR).
In particular, a compelling feature of the MGR framework is that it allows for rapid exploration of different solver design choices---notably the order in which the unknowns are reduced, different strategies to construct coarse grid operators, and a wide gamut of smoother selections.
The representative applications presented here provide motivating use cases for MGR and demonstrate its applicability to solve realistic, large-scale problems.   It also leads to portable and modular software that can be readily incorporated into existing simulators.
%
%

The paper is organized as follows. 
%
%
\Cref{sec:mgr_overview} gives a brief development of MGR and summarizes its key algorithmic components.
\Cref{sec:linear_systems} explains in detail how MGR is applied to different linear systems arising in (1) hybrid mimetic finite difference discretization of single-phase flow, (2) compositional flow with wells, and (3) hydraulic fracturing. 
In \Cref{sec:results}, we report numerical results, including weak and strong scaling performance of MGR for these applications. 
We end the paper with some concluding remarks and potential future work in \Cref{sec:conclusions}.

%% file: Section3_MGR.tex
\section{Multigrid Reduction Framework} \label{sec:mgr_overview}
The idea of MGR has had a long history, tracing back to the work of Ries and Trottenberg \cite{Ries79,Ries83}. 
Originally developed as a ``nonstandard'' geometric multigrid method (GMG), MGR introduces the use of \emph{intermediate} grids.
This gives MGR more flexibility than a GMG method, since one has the freedom to choose these intermediate grids, instead of having to rely on the typical coarsening of GMG for uniform structured grids, i.e. the grid size in each $\{x,y,z\}$-direction is decreased by a factor of 2.
For example, in \cite{Ries83}, the authors use a grid size of $\sqrt{2}h$ in between the subsequent grids.

Recently, MGR has seen a resurgence through the work on multigrid reduction in time by Falgout et al. \cite{Falgout14,Falgout17}.
It has also been studied in the context of AMG \cite{MacLachlan06,Brannick10} and other notable reduction-based multigrid methods include a nonsymmetric AMG solver developed in \cite{Manteuffel19} and an adaptive reduction approach in \cite{MacLachlan06}.
In these developments, MGR is reinterpreted as a purely algebraic method, and unlike geometric multigrid, it can be applied to general geometries and grid types.
Because of this major advantage, MGR has also been applied successfully for problems in reservoir simulation and multiphase flow in porous media with phase transitions \cite{Bui18,Wang17}.
It is also applicable to multiphysics problems such as coupled multiphase flow and geomechanics \cite{Bui20}.
The rest of this section reviews the key algorithmic components of MGR pertaining to the applications presented in this paper. 
For a more comprehensive description of the MGR framework, we refer the reader to \cite{Bui20}.

The first key idea of MGR comes from a partitioning of a matrix $A$ of size $N\times N$ into $N_C$ C-points and $N_F$ F-points, such that
\begin{linenomath}
\begin{align}
A &= \begin{pmatrix}
A_{FF} & A_{FC} \\
A_{CF} & A_{CC}
\end{pmatrix} 
= \begin{pmatrix}
I & 0 \\
A_{CF}A_{FF}^{-1} & I
\end{pmatrix} \begin{pmatrix}
A_{FF} & 0 \\
0 & S
\end{pmatrix} \begin{pmatrix}
I & A_{FF}^{-1}A_{FC}  \\
0& I
\end{pmatrix}, \label{eq:ldu_factorization}
\end{align}
\end{linenomath}
where $I$ is the identity matrix and $S = A_{CC} - A_{CF} A_{FF}^{-1} A_{FC}$ is the Schur complement.
Here, the C-points play a role analogous to the points on a coarse grid, and the F-points belong to the set that is the complement of the C-points, i.e. $N = N_C + N_F$.
Note that this partitioning is different from the one normally used in standard multigrid methods, in which the F-points correspond to all points on the fine grid, i.e. the set of F-points contains the set of C-points.
In MGR, the C-points and F-points belong to non-overlapping sets.

The second algorithmic feature that sets MGR apart from standard multigrid methods lies in the definition of the interpolation and restriction operators.
While standard multigrid methods, which specifically target elliptic problems, construct these operators based on knowledge of the near null-space, MGR derives them from the LDU factorization in \Cref{eq:ldu_factorization}. 
Ideal MGR interpolation and restriction operators can be defined as
\begin{linenomath}
\begin{align}
P = \begin{pmatrix}
-A_{FF}^{-1} A_{FC}\\
I
\end{pmatrix}, \hspace{5mm} R = \begin{pmatrix}
-A_{CF}A_{FF}^{-1} &I
\end{pmatrix}. \label{eq:ideal_iterp_restrict}
\end{align}
\end{linenomath}
Using this definition for $R$ and $P$, the coarse grid can be computed using the Petrov-Galerkin triple matrix product $S = RAP$. 
It is clear that the formulation of the ideal interpolation and restriction operators in \eqref{eq:ideal_iterp_restrict} facilitates an algebraic method, as they only require the block $A_{FF}$ to be invertible. 
From a practical standpoint, this is computationally expensive and results in a dense $P$ and $R$. 
As a result, the coarse grid is also dense and sparse approximations to $A_{FF}^{-1}$ are needed. 
One strategy that is frequently effective is to use a Jacobi approach for interpolation and an injection operator for restriction, namely
\begin{linenomath}
\begin{align}
\widetilde{P} = \begin{pmatrix}
-D_{FF}^{-1} A_{FC}\\
I
\end{pmatrix}, \hspace{5mm} 
\widetilde{R} = \begin{pmatrix}
0 &I
\end{pmatrix}. \label{eq:approximate_jacP_injR}
\end{align}
\end{linenomath}
where $D_{FF} = \text{diag}(A_{FF})$.
We note that injection is not a good option for restriction operator in traditional multigrid settings, \Red{in which the global matrix $A$ is SPD}.
\Red{An option for SPD matrices using MGR is}
\begin{linenomath}
\begin{align}
 \widetilde{R} = \begin{pmatrix}
-A_{CF}D_{FF}^{-1} & I
\end{pmatrix}, \hspace{5mm}
\widetilde{P} = \widetilde{R}^T = \begin{pmatrix}
-D_{FF}^{-1} A_{FC}\\
I
\end{pmatrix}. \label{eq:approximate_jacP_jacR}
\end{align}
\end{linenomath}
However, \Red{For multiphysics problems, the global matrix $A$ is almost never SPD. Furthermore, in these cases, the CF splitting in MGR does not need to strictly correspond to low and high frequency errors.}
Instead, as will be seen in the next section, the C and F points will be associated with physical fields, \Red{and a non-symmetric approximation to $P$ and $R$ as in \cref{eq:approximate_jacP_injR} is preferable.} 
Other options for approximating restriction and interpolation operators are summarized in \cite{Bui20}.

\begin{algorithm}
  \caption{General multi-level MGR preconditioner, $z = M_{\ell,MGR}^{-1} v$.}\label{algo:multi_level_mgr_NC}
  \begin{algorithmic}[1]
  \Function{\tt{V\_MGR}}{$A_{\ell}, v_{\ell}$}
    \If{$\ell$ is the coarsest level}
      \State $A_{\ell} z_{\ell} = v_{\ell}$ \Comment{Solve coarse-grid error problem}
    \Else
      \State $z_{\ell} = M_{\ell,glo}^{-1} v_{\ell} $ \Comment{Global Relaxation}
      \State $z_{\ell} \leftarrow z_{\ell} + Q_{\ell}^{T}M_{\ell,FF}^{-1}Q_{\ell} (v_{\ell} - Az_{\ell})$ \Comment{F-Relaxation}
      \State $r_{\ell+1} = \widetilde{R}_{\ell} (v_{\ell} - A_{\ell}z_{\ell})$ \Comment{Restrict residual}
      \State  $e_{\ell+1} = {\tt{V\_MGR}}(A_{\ell+1}, r_{\ell+1})$  \Comment{Recursion} 
      \State $e_{\ell} = \widetilde{P}_{\ell} e_{\ell+1}$ \Comment{Interpolate coarse error approximation}
      \State $z_{\ell} \leftarrow z_{\ell} + e_{\ell}$ \Comment{Apply correction}
   	\EndIf
   	\State \Return{$z_{\ell}$}
  \EndFunction
  \end{algorithmic}
\end{algorithm}

Finally, as it is a multigrid method, MGR consists of complementary processes: global smoothing, F-relaxation, and coarse-grid correction. 
Assuming the hierarchy of coarse grid operators $A_{\ell+1} = \widetilde{R}_\ell A_\ell \widetilde{P}_\ell$ has been computed at every level $\ell$, the general multi-level MGR method is shown in \Cref{algo:multi_level_mgr_NC}, where $M_{\ell,glo}^{-1}$ is the global smoother, $M_{\ell,FF}^{-1}$ is the F-relaxation smoother, and $Q_\ell = \begin{pmatrix} I & 0 \end{pmatrix}$ is the F-injection operator.
Overall, all three algorithmic components, i.e. the coarsening process, the approximations to the interpolation and restriction operators, and the choices of solvers are critical to the success of MGR.
\par \Red{In principle, as a framework, MGR provides a mechanism to quickly experiment with different combinations of these components. As such, it is possible to apply MGR to a new problem in a naive way by simply enumerating through all the possible combinations. There may be cases in which this is the only feasible approach, but  it can be non-trivial for systems with a large number of blocks. A much more efficient use of MGR is typically driven by theory or guided by physics when such knowledge is available.}
\Red{For example, the choice of solvers for the C and F blocks involved should be compatible with their properties.}
\Red{For elliptic or near-elliptic blocks, it would make sense to use an AMG V-cycle, while for mass matrices and ``hyperbolic'' blocks, cheap smoothers such as classical relaxation schemes (e.g. Jacobi, Gauss-Seidel), or variants of incomplete factorization (ILU) should be the methods of choice; and for special cases with small block sizes or algebraic constraints, a direct solver is often used.}
\par \Red{While incorporating theory and physics could help reduce the space of choices significantly, it is possible that several viable reduction strategies exist, and experimentation is still required to find the optimal configuration. This hybrid approach is the one that we take in this paper.}
%
%
%
In the next sections, we describe how \Red{a physics-guided} MGR can be used to produce efficient preconditioners for a broad range of applications of interest.

%% file: Section4_Linear_Systems.tex
\section{Representative Applications} \label{sec:linear_systems}
In this section, we describe discrete linear systems arising in three major applications considered in this work:
(1) a hybrid mimetic finite difference discretization of single-phase flow,
(2) \Red{a standard finite volume formulation of} compositional multiphase flow with wells, and
(3) \Red{a mixed finite element/finite volume scheme for} hydraulic fracturing.
For details regarding the governing equations that lead to these systems, see \ref{app:hybrid_mfvm}, \ref{app:compositional_flow}, and \ref{app:hydraulic_fracturing}.

\subsection{Hybrid MFD Method for Single-phase Flow}\label{subsec:hybrid_mfvm}
\Blue{We consider a model representing compressible single-phase flow in heterogeneous porous media.}
The (near-elliptic) PDE written in mixed form is discretized fully implicitly using a hybrid MFD scheme \cite{DaVeiga14, Lie19} involving cell-centered degrees of freedom---the cell-centered pressure, $p$---as well as face-centered degrees of freedom---the one-sided face flux, $w$, and the Lagrange multiplier, $\pi$---as described in detail in \ref{app:hybrid_mfvm}.
A key feature of the MFD scheme is that it can provide a consistent and accurate discretization for complex meshes made of highly distorted cells often encountered in flow simulation in complex geological formations \cite{Lie19,Coo_etal20}.
The scheme yields large saddle-point block linear systems \cite{Chavent86,Brezzi12} that cannot be easily handled by black-box preconditioners.
The standard linear solution strategy for these systems starts with an exact Schur complement reduction performed during assembly to eliminate the flux degrees of freedom, referred to as static condensation.
Then, scalable solvers can be obtained by applying a physics-based block preconditioner to the reduced system \cite{Frigo20,Borio20}.
As illustrated next, MGR provides a convenient framework to implement such a preconditioner in a purely algebraic fashion.

\Blue{In the compressible case, the model is nonlinear, and yields a non-symmetric hybrid MFD linear system at each nonlinear iteration.}
\Blue{It can be written in the form:}
%
\begin{linenomath}
\begin{align}
    A = \begin{pmatrix}
    A_{ww} & A_{wp} & -A_{w\pi} \\
    A_{pw} & A_{pp} & 0 \\
    A^T_{w\pi} & 0 & 0 \\
    \end{pmatrix}. \label{eq:mfd_saddle_point_system}
\end{align}
\end{linenomath}
%
%
Here,
\begin{itemize}
    \item $A_{ww}$ is a block diagonal SPD matrix. It has one diagonal block per mesh cell, and the size of each block is equal to the number of faces in the corresponding cell. On arbitrary polyhedral meshes, the blocks of $A_{ww}$ can therefore have different sizes, but for the tetrahedral meshes used in Section \ref{sec:results}, all the blocks are of size $4 \times 4$;
    \item $A_{pp}$ has the sparsity pattern of a discrete cell-centered Laplacian containing accumulation and flux term derivatives. Given that an upwinding scheme is used in the flux approximation, $A_{pp}$ is not symmetric;
    \item $A_{w \pi}$ is a rectangular matrix in which the number of rows is the number of one-sided faces---i.e., the number of cells multiplied by the number of faces in each cell---and the number of columns is the number of (unique) faces in the mesh. An interior face is associated with two one-sided faces (one for each neighboring cell) while a boundary face has only one one-sided face. In $A_{w \pi}$, there is a non-zero entry equal to -1  at the $i$-th row and $j$-th column if one-sided face $i$ belongs to face $j$; 
    \item $A_{wp}$ and $A^T_{pw}$ are rectangular matrices sharing the same sparsity pattern, with a non-zero entry at the $i$-th row and $j$-th column if one-sided face $i$ belongs to cell $j$. We note that, in the compressible case, we have $-A_{wp} \neq A^T_{pw}$.
\end{itemize}
The detailed description of each block can be found in \ref{app:hybrid_mfvm}.
\Red{We mention here that the methodology employed to compute the blocks of $A_{ww}$ determines some key properties of the discretization, such as consistency and convergence upon spatial refinement.
In this work, we consider two approaches to illustrate how the knowledge of the discretization guides the construction of the MGR strategy.
In the first approach, referred to as TPFA, each block of $A_{ww}$ is a diagonal matrix whose coefficients are computed as the inverses of the half-transmissibilities in the widely used Two-Point Flux Approximation scheme \cite{Lie19}.
The TPFA approach does not yield a convergent scheme on general meshes, but remains very attractive for its efficiency and robustness.
In the second approach, referred to as consistent, each block of $A_{ww}$ is a dense matrix computed using a local consistency condition \cite{DaVeiga14}.
This approach produces a consistent numerical scheme that is convergent upon refinement on general polyhedral meshes \cite{Brezzi05}.}

Following standard practice, we apply static condensation during assembly to eliminate the flux degrees of freedom and the corresponding discrete Darcy equations.
Specifically, utilizing the fact that $A_{ww}$ is a block-diagonal SPD matrix, an exact inversion of $A_{ww}$ can be computed and the system in \Cref{eq:mfd_saddle_point_system} can be reduced to the following $2\times 2$ block system:
\begin{linenomath}
\begin{align}
    \Bar{A} = \begin{pmatrix}
    \Bar{A}_{pp} & \Bar{A}_{p\pi} \\
    \Bar{A}_{\pi p} & \Bar{A}_{\pi\pi}
    \end{pmatrix} = 
    \begin{pmatrix}
    A_{pp} - A_{pw}A^{-1}_{ww}A_{wp} & A_{pw}A^{-1}_{ww}A_{w\pi} \\
    -A^T_{w \pi}A^{-1}_{ww}A_{wp} & A^T_{w \pi}A^{-1}_{ww}A_{w\pi}
    \end{pmatrix}. \label{eq:mfd_reduced}
\end{align}
\end{linenomath}
After static condensation, the sparsity pattern of $\Bar{A}_{pp}$ remains unchanged because $A_{pw}A^{-1}_{ww}A_{wp}$ is a diagonal matrix.
\Red{The design of an MGR strategy for this system is guided by the knowledge of the discretization, as explained next.}

\Red{In the TPFA approach, the discretization scheme guarantees that $\Bar{A}_{\pi\pi} = A^T_{w \pi}A^{-1}_{ww}A_{w\pi}$ is a diagonal matrix.
To take advantage of this, we choose the cell-centered pressure $p$ as the C-point, and the Lagrange multiplier $\pi$ as the F-point.
The diagonal structure of $\Bar{A}_{\pi\pi}$ motivates the application of a Jacobi smoother in the F-relaxation.
We use the ideal MGR operator for interpolation (\Cref{eq:ideal_iterp_restrict}) and the injection operator for restriction (\Cref{eq:approximate_jacP_injR}), yielding an exact Schur complement that is inexpensive to compute:}
\begin{linenomath}
\begin{equation}
S = \Bar{A}_{pp} - \Bar{A}_{p\pi} \Bar{A}_{\pi\pi} \Bar{A}_{\pi p}. \label{eq:schur_pp}
\end{equation}
\end{linenomath}
\Red{This reduction relies on a strong theoretical foundation since it reproduces the operations involved in the formal derivation of the cell-centered TPFA scheme \cite{Lie19}.
Despite the presence of upwinded densities in $\Bar{A}_{pp}$ in the compressible case, the Schur complement system of \Cref{eq:schur_pp} remains close to SPD and is solved with an AMG V-cycle.
Alternatively, one can solve the Schur complement system with an auxiliary space preconditioner as is done in \cite{Batista20}.
}

\Red{In the consistent approach, $\Bar{A}_{\pi\pi}$ is no longer a diagonal matrix.
As a result, it becomes challenging to approximate the Schur complement of \Cref{eq:schur_pp} without sacrificing the efficiency of the preconditioner.
Our numerical tests show that using a diagonal approximation of $\Bar{A}_{\pi\pi}$ to approximate $S$ does not produce satisfactory results.
This motivates the use of another reduction order in MGR.}
Following the Schur complement approach in \cite{Lie19,Frigo20,Borio20}, we choose the Lagrange multiplier $\pi$ as the C-point, and the cell-centered pressure $p$ as the F-point.
For the F-relaxation, we \Red{still} use a simple Jacobi smoother as we find it effective for the non-symmetric upwinded flux approximation in $A_{pp}$.
Both the interpolation and restriction operators are \Red{now} chosen as in \Cref{eq:approximate_jacP_injR}, yielding the following approximate Schur complement:
\begin{linenomath}
\begin{equation}
S = \Bar{A}_{\pi \pi} - \Bar{A}_{\pi p} \text{diag}( \Bar{A}_{pp} )^{-1} \Bar{A}_{p \pi}. \label{eq:schur_pipi}
\end{equation}
\end{linenomath}
This approximate Schur complement system matrix is not symmetric but remains close to SPD.
The Schur complement system is solved with an AMG V-cycle.
We refer the reader to \cite{Dobrev19}, where the use of AMG to a similarly reduced system is motivated.
\Red{From now on, we denote MGR\_X with X=\{P,$\Pi$\} as the reduction strategy that assigns variable X as the F-points.
We note that MGR\_$\Pi$ can also be used to precondition the linear systems resulting from the TPFA approach, but this requires solving a larger (face-based) Schur complement system and is therefore more expensive that MGR\_P (see Section~\ref{sec:results}).}
  
%

\subsection{Compositional Models with Wells} \label{subsec:compositional_flow}

In this section, we consider two-phase, two-component compositional flow in the context of CO$_2$ injection for permanent storage in heterogeneous geological formations.
The complete formulation, described in \ref{app:compositional_flow}, leads to a highly nonlinear problem governing both the near-elliptic evolution of the pressure field and near-hyperbolic evolution of the CO$_2$ plume.
The reservoir flow problem is coupled fully implicitly to a multi-segment well flow model, governed by a hyperbolic PDE and controlled by a set of algebraic constraints.
The system considered here therefore involves complex physical processes operating on different time scales, as the propagation of the CO$_2$ front in the reservoir is often relatively slow compared to fast-rate flow in the wellbore and its immediate surroundings.
It also contains both discretized PDE equations and algebraic constraint equations.
For this problem, we employ a standard finite volume discretization with the backward Euler method, and linearize the discrete nonlinear system using Newton's method.
%

The most common strategy to scalably precondition this problem relies on the \textit{Constrained Pressure Residual} (CPR) multistage preconditioning technique.
Originally developed in \cite{Wallis83, Wallis85}, this approach has been extended and applied to a wide range of problems \cite{Cao05,Gries14,Lacroix03,Liu15,Scheichl03,Stueben07,Zhou12}.  
This family of preconditioners can be considered as a variant of physics-based block preconditioners, as they exploit knowledge of the physics to break the coupled linear system into smaller subproblems that can be solved separately by standard techniques.
%
%
%
Here, we demonstrate how one can apply MGR to handle this system.
\Red{We note that the MGR reduction strategy we employ shares some similarity to the established CPR preconditioners, but it also has some distinctive features.}

The linear system under consideration takes the following block form:
%
%
\begin{linenomath}
\begin{align}
A = 
\begin{pmatrix}
A_{rr} & A_{rw} \\
A_{wr} & A_{ww}
\end{pmatrix} \;, \label{eq:general_compositional_system}
\end{align}
\end{linenomath}
where the diagonal blocks $A_{rr}, A_{ww}$ represent compositional flow in the reservoir and the wells, respectively, and the off-diagonal blocks are the coupling between the reservoir and the wells.
Block $A_{rr}$ itself is a block matrix,
\begin{linenomath}
\begin{align}
A_{rr} = 
\begin{pmatrix}
A^r_{pp} & A^r_{p\rho_1} & A^r_{p\rho_2} \\
A^r_{\rho_1 p} & A^r_{\rho_1 \rho_1} & A^r_{\rho_1 \rho_2} \\
A^r_{\rho_2 p} & A^r_{\rho_2 \rho_1} & A^r_{\rho_2 \rho_2} 
\end{pmatrix} \;, \label{eq:compositional_reservoir_block}
\end{align}
\end{linenomath}
and it corresponds to three physical variables, i.e. reservoir pressure $p^r$ and reservoir component densities $\rho_1^r$, and $\rho_2^r$.
A component density of component $c$ (either CO$_2$ or H$_2$O), $\rho_c$, is defined as the mass of component $c$ per unit volume of fluid mixture.
The well block, $A_{ww}$, is block-diagonal since we assume that wells do not interact with each other.
Given that we model the flow in the wellbore using a PDE, each diagonal block of $A_{ww}$ has a similar structure to $A_{rr}$, albeit with an additional equation for the control of the pumping rate, leading to a $4\times 4$ block matrix that corresponds to four variables $p^w, \rho_1^w, \rho_2^w, \sigma^w$, i.e. well pressure, well component densities, and well rate between segments.
%
Even though these are four separate physical quantities, we group them together and treat them as a single variable $w = \{p^w, \rho_1^w, \rho_2^w, \sigma^w\}$ in this work.
Thus, in the context of MGR, it is most convenient that we view the system in \Cref{eq:general_compositional_system} as
\begin{linenomath}
\begin{align}
A = 
\begin{pmatrix}
A_{pp} & A_{p\rho_1} & A_{p\rho_2} & A_{pw} \\
A_{\rho_1 p} & A_{\rho_1 \rho_1} & A_{\rho_1 \rho_2} & A_{\rho_1 w} \\
A_{\rho_2 p} & A_{\rho_2 \rho_1} & A_{\rho_2 \rho_2} & A_{\rho 2 w} \\
A_{wp} & A_{w\rho_1} & A_{w\rho_2} & A_{ww}
\end{pmatrix} \;. \label{eq:compositional_system}
\end{align}
\end{linenomath}
Note that since we no longer need to distinguish between the reservoir and well pressures and densities, the superscript $r$ is dropped.
To motivate the design of an MGR strategy, we emphasize the following properties of $A$:
\begin{itemize}
    \item $A_{pp}$ is the pressure block that has the structure of a discrete elliptic operator;
    \item $A_{\rho_1\rho_1}$ has the structure of a discrete time-dependent hyperbolic problem;
    \item $A_{\rho_2 p}$, $A_{\rho_2\rho_1}$, $A_{\rho_2\rho_2}$ are diagonal matrices describing the cell-based volume constraint equations used in the compositional formulation to close the system of PDEs;
    \item $A_{ww}$ is the well block containing both discretized PDEs and constraint equations for pressure and rate controls. Such algebraic constraints can lead to the existence of zeros on the diagonal of this block.
\end{itemize}
\Red{Some authors such as \cite{Cao15} choose to eliminate the well block $A_{ww}$ first by inverting it exactly and then solving the reduced Schur complement system containing only variables in the reservoir with CPR.}
\Red{Here, we employ a 4-level MGR strategy that keeps the well equations until the coarsest grid.}
For all the levels, the interpolation and restriction operators specified in \Cref{eq:approximate_jacP_injR} are used.
\Red{Regarding the reduction order, we choose the reservoir densities $\rho_2$ and as F-points for the first level. As noted above, as these densities correspond to the volume constraint equations and the blocks $A_{\rho_2 p}$, $A_{\rho_2\rho_1}$, $A_{\rho_2\rho_2}$ are diagonal, the MGR reduction for this first level is exact with the choice of interpolation and restrictions in \Cref{eq:approximate_jacP_injR}. The F-relaxation should also be exact with a single Jacobi iteration.}
\Red{For the second level, we eliminate the reservoir densities $\rho_1$ that correspond to the advective component of the system. As a result, Jacobi is also an appropriate choice as a F-relaxation smoother.}
%
%
For the third level, we eliminate the reservoir pressure. This means that the $A_{FF}$ matrix is a perturbation of the elliptic pressure block $A_{pp}$.
Thus, an AMG V-cycle is applied for F-relaxation.
\Red{Note that due to the use of interpolation and restrictions operator in \Cref{eq:approximate_jacP_injR}, the $A_{FF}$ matrix is exactly the pressure system in a CPR approach and this reduction step can be seen as solving the Schur complement pressure system.}
After this reduction step, the global system is reduced onto a much smaller coarse-grid system that corresponds to the well equations.
Because of the small size of the well block, we use a direct method with SuperLU\_Dist \cite{Li03} to solve the coarse grid.

One distinctive feature of the system in \Cref{eq:compositional_system} is the presence of zeros on the diagonal of $A_{ww}$ due to the constraints imposed on the fluid pressure and pumping rate for each well.
The fact that the system contains both discretized PDEs and algebraic constraints makes it very challenging to solve, even for state-of-the-art methods like CPR.
CPR strategies usually involve a global smoothing step following a pressure correction solve, and their effectiveness relies significantly on the convergence of these two processes.
The most popular global smoother for CPR approaches is based on ILU.
The presence of zeros on the diagonal could cause stability issues for ILU, while other methods such as Jacobi and Gauss-Seidel cannot be applied at all.
A notable advantage of our approach is that it does not require a global relaxation step.
The separation of the reservoir and well variables by the C-F splitting ensures that the F-relaxation is applied only to reservoir variables and avoids the zeros on the diagonal in the well equations.

\subsection{Hydraulic Fracturing} \label{subsec:hydraulic_fracturing}
The last example focuses on the equilibrium solve at the heart of a hydraulic fracturing simulation \cite{Adachi:2007,Settgast:2017}.  The governing equations (\ref{app:hydraulic_fracturing}) are discretized in space using a finite element scheme for the momentum balance equation modeling the deformation of the porous medium, and a standard finite volume scheme for the mass conservation equation of the fracturing fluid.  In this example, no leakoff of fluids from the fracture to the surrounding rock is considered.
Time integration is handled using a backward Euler scheme.
The resulting linearized Jacobian system has the natural block structure
\begin{linenomath}
\begin{equation}
A = \begin{pmatrix}
A_{uu} & A_{up} \\
A_{pu} & A_{pp} 
\end{pmatrix} \,. \label{eq:hydraulic_fracturing_linear_system}
\end{equation}
\end{linenomath}
Here, $A_{uu}$ is an elastic operator from the volumetric mesh and $A_{pp}$ is a Laplacian-like operator from the lower-dimensional fluid mesh.
Both are discrete elliptic operators.
The off-diagonal $A_{up}$ captures coupling between the pressure exerted on the fracture walls and the surrounding rock deformation, while $A_{pu}$ captures the impact of the changing fracture aperture on the fluid storage term of the mass balance equation.
%
%
%

For coupled fracture mechanics and fracture fluid flow, very few large-scale solvers exist.
In fact, in many fully-implicit three-dimensional geomechanical fracture simulators, including the one developed in \cite{Zheng19} and the one used in  \cite{Settgast:2017}, distributed direct solvers \cite{Sala06,Li03} are widely used.
Regarding iterative methods, preconditioners have been proposed for alternative discretization approaches to the underlying governing equations. 
For example, in \cite{Peirce05,Peirce06}, the same author considered a specialized geometric multigrid method and an incomplete factorization (ILU) approach for a localized Jacobian operator resulting from a coupled integro-partial differential equations system that is discretized using a combination of collocation and finite volume methods.

To the best of our knowledge, no physics-based preconditioners have been developed, and none of the current solvers is capable of tackling the problem of interest specified in \Cref{eq:hydraulic_fracturing_linear_system}.
Thus, in this context, MGR plays two simultaneous roles: (1) as a high-performance algebraic preconditioner, and (2) as a framework to help guide future development of specialized physics-based counterparts.
Due to the natural block structure of the problem, we propose a two-level MGR method with a C-F splitting that aligns with the physical variables.
One can choose two possible reduction strategies: reducing the mechanics block $A_{uu}$ on to the pressure block $A_{pp}$, or vice versa.
\Blue{For the choices of approximate interpolation and restriction operators, either \Cref{eq:approximate_jacP_injR},  \Cref{eq:approximate_jacP_jacR}, or approximate inverse can be used. However, in our experiments, we found that the first approach obtains the best convergence and produces the sparsest coarse grid systems at the same time. Regarding the choices of F-relaxation smoother and coarse-grid solver, the fact that both the diagonal blocks are derived from an elliptic operator helps limit the choices. We use variants of an AMG V-cycle, as detailed below.}
In \Cref{sec:results}, we present the results for both strategies.
Here, we illustrate the first reduction strategy in detail and discuss the differences and advantages of each approach.

Reducing the mechanics block $A_{uu}$ first implies $A_{FF} \equiv A_{uu}$.
Using the interpolation and restriction operators in \Cref{eq:approximate_jacP_injR}, the coarse grid can be computed as
\begin{linenomath}
\begin{align}
    \widetilde{S} = \widetilde{R}A\widetilde{P} = A_{pp} - A_{pu}D_{uu}^{-1}A_{up},
\end{align}
\end{linenomath}
in which $D_{uu} = \text{diag}(A_{uu})$.
For the F-relaxation step, since the $A_{uu}$ block is a discrete elliptic operator coming from discretizing the elasticity equation, we use one AMG V-cycle.
One can show that a block diagonal approximation, in which there is no coupling between displacement in each direction $\{x,y,z\}$ (also called \textit{separate displacement component} \cite{Axelsson78}), is spectrally equivalent to the original $A_{uu}$ \cite{Blaheta94,Gustafsson98}.
Thus, we employ the unknown-based approach of AMG \cite{Stueben01} for the three variables ($u_x, u_y$ and $u_z$).
The unknown approach is essentially a block diagonal preconditioner, in which each variable $u_x, u_y, u_z$ has its own AMG hierarchy.
The difference is that in the unknown approach, the coarse grid still includes the off-diagonal coupling blocks between the variables.
Regarding the coarse grid operator, since it is also a perturbation of a discrete elliptic operator, i.e. $A_{pp}$, we apply one AMG V-cycle for the coarse-grid correction step.
The advantage of this splitting is that the coarse grid, associated with the fluid pressure in the fractures, is significantly smaller than the global problem.
Furthermore, the larger mechanics block solved in the F-relaxation step is unperturbed, making it ideal for AMG.
Additionally, unless the fractures propagate, the mechanics block is kept fixed for each Newton solve, and using this reduction scheme allows for reuse of its AMG hierarchy, thus reducing the setup cost.


However, it may make more sense to solve for the pressure first in the F-relaxation step, as the flow in the fractures are driven by the fluid pressure.
Thus, letting $A_{FF} \equiv A_{pp}$ would lead to a coarse grid that is more representative of the physics.
Nevertheless, the resulting coarse grid, which incorporates the coupling from the reduction, may be more difficult for AMG to handle.
Another disadvantage of this reduction scheme stems from the fact that the pressure matrix representing fluid flow is updated at every Newton iteration.
The coarse grid, which is a perturbation of the mechanics block $A_{uu}$ in this case, also changes at every Newton iteration and needs to be re-setup as a result.
Similar to the first reduction approach, we use one AMG V-cycle for the F-relaxation, and the unknown-based AMG approach with three variables for the coarse grid correction.
From now on, we denote MGR\_X with X=\{U,P\} as the reduction strategy that assigns variable X as the F-points. In the next section, we present and discuss the results for both approaches.

%% file: Section5_Results.tex
\section{Numerical Results}\label{sec:results}
\subsection{Implementation Details}
MGR is implemented as a separate solver in \textit{hypre} - a high performance linear solver and preconditioner package \cite{Falgout02}. 
As such, MGR can leverage other solvers and infrastructure within \textit{hypre}, such as BoomerAMG \cite{Henson00} and parallel ILU solvers. 
All the numerical experiments were run on Quartz, a cluster at the Lawrence Livermore Computing Center with 1344 nodes containing two Intel Xeon E5-2695 18-core processors sharing 128 GB of memory on each node with Intel Omni-Path interconnects between nodes.
We use pure MPI-based parallelism.

Unless specified otherwise, a default configuration of BoomerAMG is used to construct an AMG V-cycle. 
Specifically, at each level, the smoother consists of one sweep of hybrid forward $L_1$-Gauss-Seidel \cite{Baker11} for the down cycle and one sweep of hybrid backward $L_1$-Gauss-Seidel for the up cycle. 
The coarsest grid is solved directly with Gaussian elimination. 
For a scalar problem, i.e. pressure solves, a Hybrid Modified Independent Set (HMIS) coarsening strategy \cite{DeSterck06} is used. 
For the system case, i.e. the mechanics block, an unknown approach is used with one level of aggressive coarsening. No global smoother is needed for all the block systems targeted here.
The MGR preconditioner is always used in conjunction with right-preconditioned GMRES.

The different physical formulations used for the test applications are all implemented in GEOSX\footnote{The public GitHub repository can be found at \url{https://github.com/GEOSX/GEOSX}.} \cite{geosx}.

\subsection{Hybrid Mimetic Finite Difference for Single-Phase Flow}\label{subsec:hybrid_mfvm_results}

\begin{figure}
  \small
  \centering
  \begin{subfigure}[b]{0.475\textwidth}
    \centering
    \includegraphics[scale=1]{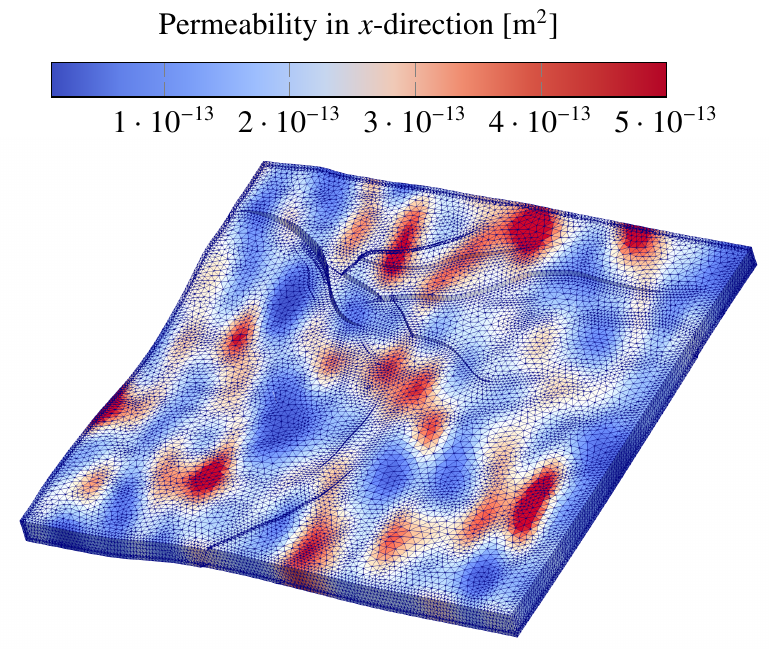}
    \caption{}
    \label{fig:HI24L_porosity_permeability}
  \end{subfigure}
  \hfill
  \begin{subfigure}[b]{0.475\textwidth}
    \centering
    \includegraphics[scale=1]{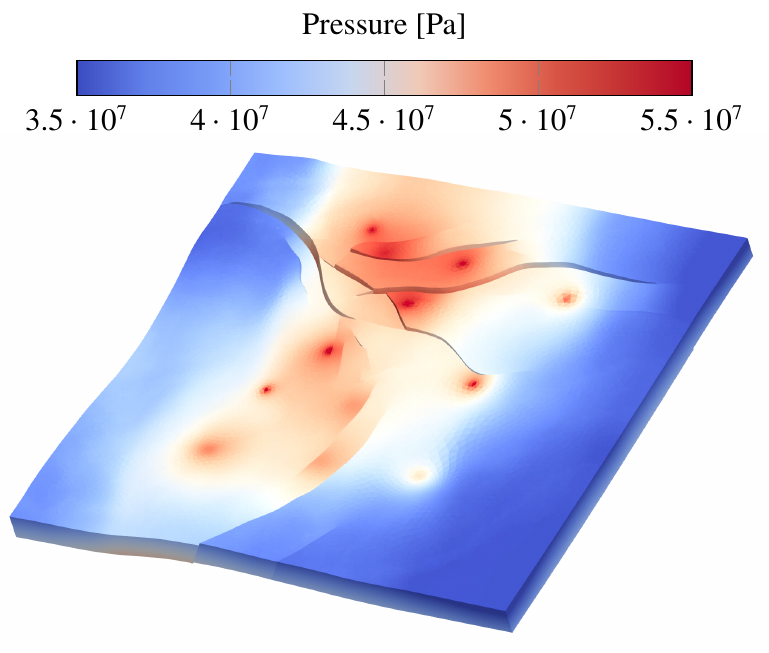}
    \caption{}
    \label{fig:HI24L_pressure_map}
  \end{subfigure}
  \caption{(a) Permeability field used in the single-phase and compositional HI24L problems. We refer the reader to \cite{Mazuyer20} for the exact location of the faults, treated as impermeable barriers to flow in this work. The porosity field, not shown here, is correlated to the permeability field. The mesh displayed in this figure is the mid-size mesh consisting of 1,560,722 tetrahedra. (b) Pressure field after 6 years of injection in the compressible single-phase problem. The single-phase problem has been set up such that the high-pressure areas correspond to the location of the twelve wells of the compositional test case.}
  \label{fig:three graphs}
\end{figure}

Here, we study the performance of MGR for a compressible single-phase simulation using two fully unstructured tetrahedral meshes modeling a domain referred to as the HI24L reservoir.
This domain represents a depleted gas field off the Gulf of Mexico considered as a candidate for industrial-scale geological CO$_2$ sequestration \cite{Deangelo19}.
It is 19-km long, 18.3-km wide, approximately 1.0-km thick.
The two meshes include a medium-sized version consisting of 1,560,722 tetrahedra, and a larger one made of 7,709,462 tetrahedra. Both meshes have 32 faults assumed to be impermeable \cite{Mazuyer20}.
On the mid-size mesh, the discretization described in \ref{app:hybrid_mfvm} involves 1,560,722 cell-centered dofs and 3,494,826 face-centered dofs, for a total of 5,055,548 dofs.
On the larger mesh, the scheme consists of 7,709,462 cell-centered dofs and 15,823,665 face-centered dofs, for a total of 23,533,127 dofs.
The heterogeneous and anisotropic permeability field is shown in \Cref{fig:HI24L_porosity_permeability}.

We set a hydrostatic initial pressure field and impose no-flow boundary conditions.
For this problem, we do not use a well model, which will be considered in the next section.
Instead, we simulate injection by imposing a fixed pressure of $6.0 \times 10^7$ Pa at 12 different locations in the domain.
We simulate $2.0 \times 10^8$ s (approximately 6 years) of pressure evolution.
\Red{In the physical model considered here, the nonlinearity only comes from the slight fluid compressibility, which makes it possible to take large time steps.}
The initial time step of $10^5$ s is doubled whenever Newton's method converges in fewer than 9 Newton iterations, and the largest time step size taken during the simulation is $10^8$ s (approximately 3 years).
The final pressure field, shown in \Cref{fig:HI24L_pressure_map}, illustrates the pressure buildup taking place in the neighborhood of the wells.

\begin{table}
    \small
    \centering
    \caption{Strong scaling for MGR\_P and MGR\_$\Pi$ for the compressible single-phase HI24L problem discretized with the TPFA hybrid MFD scheme. The test case involves 5.1M degrees of freedom and uses \textit{tol} = $10^{-7}$.}
    \begin{tabular}{lcccccccccc}
        \toprule
        Number of cores & \multicolumn{2}{c}{8} & \multicolumn{2}{c}{16} & \multicolumn{2}{c}{32} & \multicolumn{2}{c}{64} & \multicolumn{2}{c}{128} \\
        \midrule
        Avg. problem size per core & \multicolumn{2}{c}{$\sim$632,000} & \multicolumn{2}{c}{$\sim$316,000} & \multicolumn{2}{c}{$\sim$158,000} & \multicolumn{2}{c}{$\sim$79,000} & \multicolumn{2}{c}{$\sim$39,500} \\
        \cmidrule(l){2-3} \cmidrule(l){4-5} \cmidrule(l){6-7} \cmidrule(l){8-9} \cmidrule(l){10-11}
        MGR strategy & P & $\Pi$ & P & $\Pi$ & P & $\Pi$ & P & $\Pi$ & P & $\Pi$ \\
        Avg. number of iterations & 10.06 & 9.12 & 10.00 & 9.12 & 10.06 & 9.30 & 10.12 & 9.21 & 10.09 & 9.12 \\
        Avg. setup time (s) & 3.04 & 2.19 & 1.34 & 0.91 & 0.76 & 0.43 & 0.38 & 0.20 & 0.20 & 0.11 \\
        Avg. solve time (s) & 2.04 & 1.67 & 0.80 & 0.62 & 0.46 & 0.33 & 0.23 & 0.16 & 0.12 & 0.08 \\
        \bottomrule
    \end{tabular}
    \label{tab:hybrid_mfd_strong_scaling_small_tpfa}
\end{table}

\begin{table}
    \small
    \centering
    \caption{Strong scaling for MGR\_P and MGR\_$\Pi$ for the compressible single-phase HI24L problem discretized with the TPFA hybrid MFD scheme. The test case involves 23.5M degrees of freedom uses \textit{tol} = $10^{-7}$.}
    \begin{tabular}{lcccccccccc}
        \toprule
        Number of cores & \multicolumn{2}{c}{32} & \multicolumn{2}{c}{64} & \multicolumn{2}{c}{128} & \multicolumn{2}{c}{256} & \multicolumn{2}{c}{512} \\
        \midrule
        Avg. problem size per core & \multicolumn{2}{c}{$\sim$700,000} & \multicolumn{2}{c}{$\sim$350,000} & \multicolumn{2}{c}{$\sim$175,000} & \multicolumn{2}{c}{$\sim$87,500} & \multicolumn{2}{c}{$\sim$43,750} \\
        \cmidrule(l){2-3} \cmidrule(l){4-5} \cmidrule(l){6-7} \cmidrule(l){8-9} \cmidrule(l){10-11}
        MGR strategy & P & $\Pi$ & P & $\Pi$ & P & $\Pi$ & P & $\Pi$ & P & $\Pi$ \\
        Avg. number of iterations & 11.97 & 9.56 & 12.12 & 9.65 & 12.09 & 9.62 & 12.09 & 9.88 & 12.26 & 9.65 \\
        Avg. setup time (s) & 3.25 & 2.34 & 1.51 & 1.05 & 0.71 & 0.48 & 0.39 & 0.25 & 0.24 & 0.18  \\
        Avg. solve time (s) & 2.50 & 1.68 & 1.28 & 0.91 & 0.64 & 0.38 & 0.32 & 0.18 & 0.19 & 0.12 \\
        \bottomrule
    \end{tabular}
    \label{tab:hybrid_mfd_strong_scaling_large_tpfa}
\end{table}

\begin{table}
    \small
    \centering
    \caption{Strong scaling for MGR\_P for the compressible single-phase HI24L problem discretized with the consistent hybrid MFD scheme. The test case involves 5.1M degrees of freedom and uses \textit{tol} = $10^{-7}$.}
    \begin{tabular}{lccccc}
        \toprule
        Number of cores & 8  & 16 & 32 & 64 & 128 \\
        \midrule
        Avg. problem size per core & $\sim$632,000 & $\sim$316,000 & $\sim$158,000 & $\sim$79,000 & $\sim$39,500\\
        Avg. number of iterations & 12.97 & 12.82 & 12.82 & 13.42 & 12.91 \\
        Avg. setup time (s) & 1.91 & 0.89 & 0.45 & 0.20 & 0.11 \\
        Avg. solve time (s) & 2.26 & 1.05 & 0.51 & 0.27 & 0.15 \\
        \bottomrule
    \end{tabular}
    \label{tab:hybrid_mfd_strong_scaling_small_consistent}
\end{table}

\begin{table}
    \small
    \centering
    \caption{Strong scaling for MGR\_P for the compressible single-phase HI24L problem discretized with the consistent hybrid MFD scheme. The test case involves 23.5M degrees of freedom and uses \textit{tol} = $10^{-7}$.}
    \begin{tabular}{lccccc}
        \toprule
        Number of cores & 32 & 64 & 128 & 256 & 512 \\
        \midrule
        Avg. problem size per core & $\sim$700,000 & $\sim$350,000 & $\sim$175,000 & $\sim$87,500 & $\sim$43,750\\
        Avg. number of iterations & 15.58 & 15.70 & 15.55 & 15.79 & 15.71 \\
        Avg. setup time (s) & 3.34 & 1.67 & 0.90 & 0.43 & 0.22 \\
        Avg. solve time (s) & 3.62 & 1.74 & 0.82 & 0.41 & 0.24 \\
        \bottomrule
    \end{tabular}
    \label{tab:hybrid_mfd_strong_scaling_large_consistent}
\end{table}

\Red{We first consider the performance of the MGR\_P and MGR\_$\Pi$ strategies applied to the TPFA hybrid MFD approach (see Section~\ref{subsec:hybrid_mfvm}).
The strong scalability results for the two meshes are in \Cref{tab:hybrid_mfd_strong_scaling_small_tpfa,tab:hybrid_mfd_strong_scaling_large_tpfa}.
Both approaches require a small number of GMRES iterations, but MGR\_$\Pi$ is consistently more efficient than MGR\_P in terms of linear iteration count and CPU timings.
The difference is particularly significant for the last (and most difficult) time step of the fine simulation, for which MGR\_$\Pi$ requires, on average, 38\% fewer iterations than MGR\_P.
As explained in Section~\ref{subsec:hybrid_mfvm}, this is because MGR\_$\Pi$ fully exploits the structure of the system arising from the TPFA approach and only involves a cell-centered solve on the coarse grid.}

\Red{However, the reduction order in MGR\_$\Pi$ makes it more challenging to handle the structure of the system arising from the consistent hybrid MFD approach and can lead to slow linear convergence when this discretization is selected. 
Therefore, for the consistent approach, we only report the results for obtained with MGR\_P. 
\Cref{tab:hybrid_mfd_strong_scaling_small_consistent,tab:hybrid_mfd_strong_scaling_large_consistent} show that the number of GMRES iterations with MGR\_P for the consistent approach is larger than for the TPFA approach due to the higher complexity of the linear systems, but remains quite small---with about 13 and 15 iterations for the small and large mesh, respectively.}

\Red{These results confirm that the knowledge of the discretization can greatly help the elaboration of an optimal MGR strategy for a given physical problem.
To conclude this part, we note that the results of \Cref{tab:hybrid_mfd_strong_scaling_small_tpfa,tab:hybrid_mfd_strong_scaling_large_tpfa,tab:hybrid_mfd_strong_scaling_small_consistent,tab:hybrid_mfd_strong_scaling_large_consistent} indicate good algorithmic weak scaling of MGR, i.e. the number of iterations is independent of the mesh size.
In terms of strong scaling, MGR achieves near optimal performance with respect to both setup and solve time.}

\subsection{Compositional Flow with Wells}

\begin{figure}[t!]
    
    \small
    \centering
    \begin{subfigure}[b]{0.31\textwidth}
      \centering
      \includegraphics[scale=1]{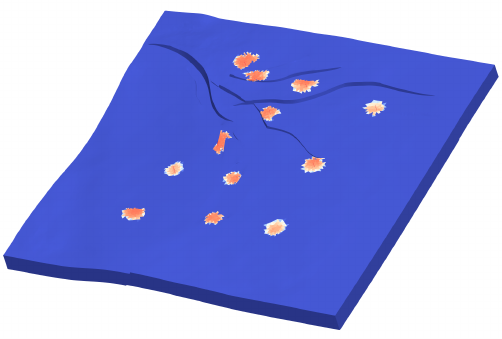}
      \caption{T $\approx$ 1 year}
    \end{subfigure}
    \hfill
    \begin{subfigure}[b]{0.31\textwidth}
      \centering
      \includegraphics[scale=1]{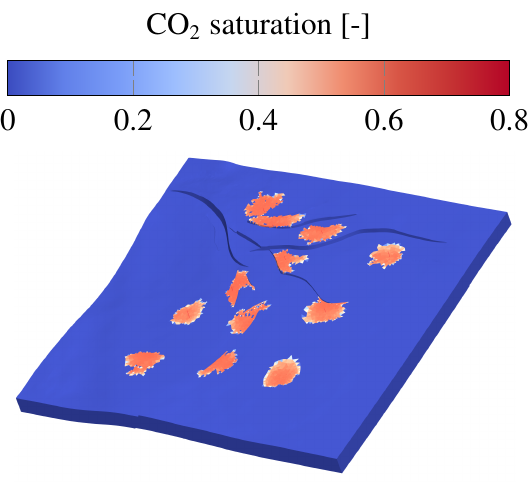}
      \caption{T $\approx$ 3 years}
    \end{subfigure} 
    \hfill
    \begin{subfigure}[b]{0.31\textwidth}
      \centering
      \includegraphics[scale=1]{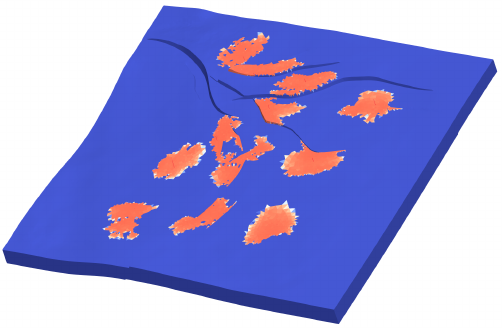}
      \caption{T = 10 years}
    \end{subfigure}    
    
    \caption{CO$_2$ saturation maps in the compositional HI24L problem, where the saturation is defined as the fraction of the pore space occupied by the CO$_2$ phase. The CO$_2$ plume is lighter than the resident brine and accumulates at the top of the domain as the simulation progresses. The shape of the plumes is constrained by the impermeable faults present in the domain.}
    \label{fig:HI24L_saturation_map}
\end{figure}

We now assess the performance of MGR for an isothermal compositional simulation of CO$_2$ injection with 12 multi-segment wells.
Again, we use the two meshes with a heterogeneous and anisotropic permeability field for the HI24L block described in the previous section. For this case, the total numbers of dofs for the two meshes are approximately 4.7 and 23.1 million, respectively.
We use a quadratic (Corey) relative permeability model, and all the other solid and fluid properties are taken from \cite{Buscheck12}.
No-flow (Neumann) boundary conditions are used to represent a closed system.
The wells are placed at the locations where pressure was imposed in the previous test case.
Each well is discretized with 10 segments and injects CO$_2$  at a rate of 10 kg/s (315,360 metric tonnes per year) through three perforations.
The compositional formulation and TPFA finite-volume discretization are reviewed in \ref{app:compositional_flow}.
We simulate 10 years of injection.
The initial time step is set to 100 s and is doubled when the Newton solver converges in fewer than 9 iterations.
Since the underlying physical model considered here is highly nonlinear, we set a maximum time step size of 20 days to avoid Newton convergence failures.
The saturation map of \Cref{fig:HI24L_saturation_map} shows that the CO$_2$ plume, initially injected at the bottom of the domain, migrates upward under the action of buoyancy forces and eventually accumulates at the top of the reservoir.

\begin{table}[!hb]
    \centering
    \small
    \caption{Strong scaling of MGR for the compositional HI24L problem with 4.7M degrees of freedom and \textit{tol} = $10^{-6}$.}

    \begin{tabular}{lccccc}
        \toprule
        Number of cores & 8 & 16 & 32 & 64 & 128 \\
        \midrule
        Avg. problem size per core & $\sim$588,000 & $\sim$294,000 & $\sim$147,000 & $\sim$73,500 & $\sim$36,800 \\
        Avg. number of iterations & 11.32 & 11.81 & 11.21 & 11.44 & 11.47 \\
        Avg. setup time per linear solve (s) & 3.80 & 1.82 & 0.80 & 0.38 & 0.19 \\
        Avg. solve time per linear solve (s) & 3.93 & 1.88 & 0.91 & 0.45 & 0.24 \\
        \bottomrule
    \end{tabular}
    \label{tab:mgr_hi24l_small}
\end{table}

\begin{table}[ht]
    \centering
    \small
    \caption{Strong scaling of MGR for the compositional HI24L problem with 23.1M degrees of freedom and \textit{tol} = $10^{-6}$.}
    \begin{tabular}{lccccc}
        \toprule
        Number of cores & 32 & 64 & 128 & 256 & 512 \\
        \midrule
        Avg. problem size per core & $\sim$700,000 & $\sim$350,000 & $\sim$175,000 & $\sim$87,500 & $\sim$43,750\\
        Avg. number of iterations & 14.85 & 14.67 & 15.20 & 15.00 & 15.25  \\
        Avg. setup time per linear solve (s) & 7.90 & 3.55 & 1.71 & 0.89 & 0.56  \\
        Avg. solve time per linear solve (s) & 5.95 & 2.89 & 1.29 & 0.71 & 0.42 \\
        \bottomrule
    \end{tabular}
    \label{tab:mgr_hi24l_large}
\end{table}

\par \Cref{tab:mgr_hi24l_small} shows the strong scaling results of MGR for the medium-sized HI24L CO$_2$/brine compositional flow problem.
Similar to the single-phase MFD case, MGR setup phase achieves almost optimal scaling, capable of reducing the run time with increasing number of cores down to as low as approximately 12,000 degrees of freedom per core.
In terms of the solve time, MGR also scales perfectly with only a minor decrease in efficiency for 128 cores.

\Cref{tab:mgr_hi24l_large} shows the strong scaling results of MGR for the large HI24L CO$_2$/brine compositional flow problem. 
Compared to the medium-sized example, the problem size increases by a factor of 5, and there is a slight increase in the number of linear iterations but the growth rate is well under control. 
Overall, the average number of linear iterations is still quite small for such a large and challenging problem. 
In terms of setup and solve time, MGR also exhibits near optimal strong scaling as observed before, except for a minor decrease in efficiency for the 512 cores.

\subsection{Hydraulic Fracturing}

\begin{figure}[b!]
    \small
    \centering
    \begin{subfigure}[b]{0.31\textwidth}
      \centering
      \includegraphics[scale=1]{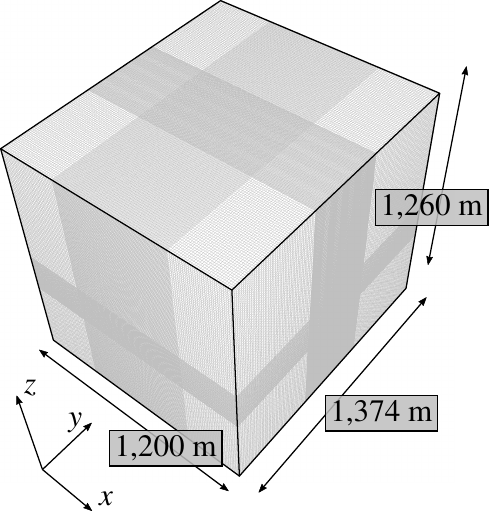}
      \caption{Computational mesh.}
      \label{fig:domain_mesh}
    \end{subfigure}
    \hfill
    \begin{subfigure}[b]{0.31\textwidth}
      \centering
      \includegraphics[scale=1]{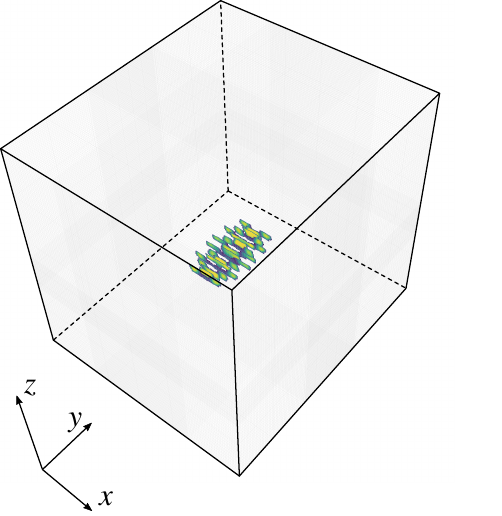}
      \caption{Full domain view.}
      \label{fig:all_fractures_aperture_full}
    \end{subfigure} 
    \hfill
    \begin{subfigure}[b]{0.31\textwidth}
      \centering
      \includegraphics[scale=1]{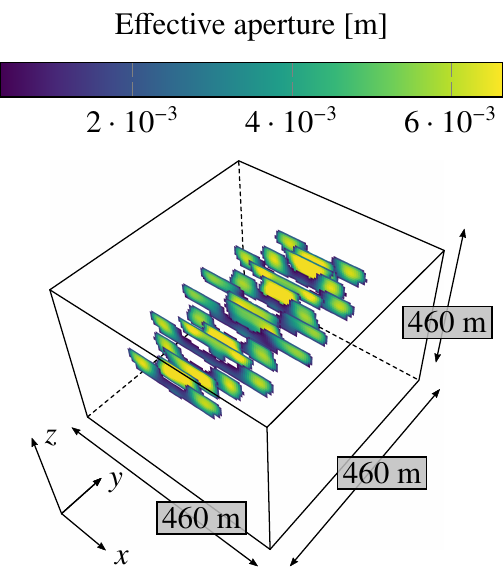}
      \caption{Zoomed in view.}      
      \label{fig:all_fractures_aperture_zoom}
    \end{subfigure} 
    \caption{Effective aperture of all the fractures at time T = 908s.}
\end{figure}

\begin{figure*}[t!]
  \small
    \small
    \centering
    \begin{subfigure}[b]{0.31\textwidth}
      \centering
      \includegraphics[scale=1]{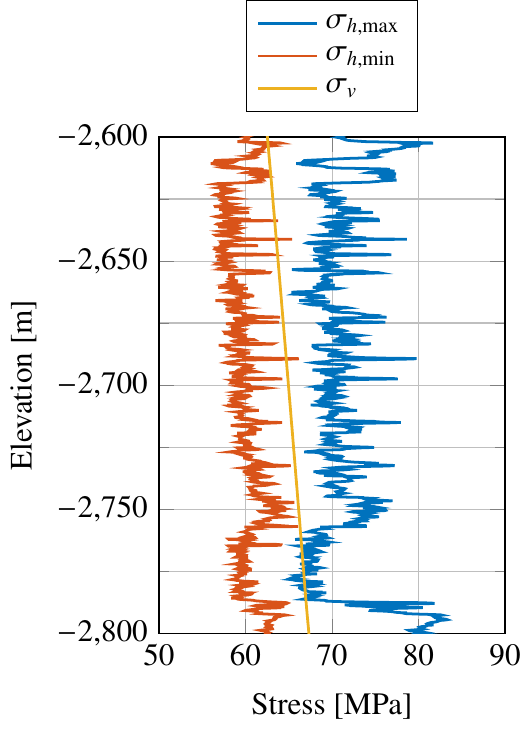}
    \end{subfigure}
    \hfill
    \begin{subfigure}[b]{0.31\textwidth}
      \centering
      \includegraphics[scale=1]{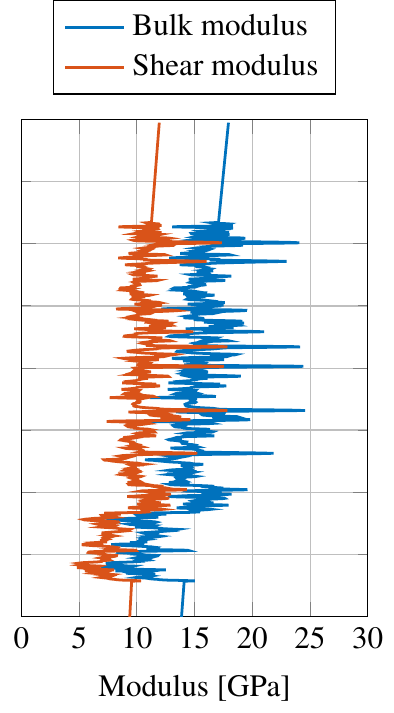}
    \end{subfigure} 
    \hfill
    \begin{subfigure}[b]{0.31\textwidth}
      \centering
      \includegraphics[scale=1]{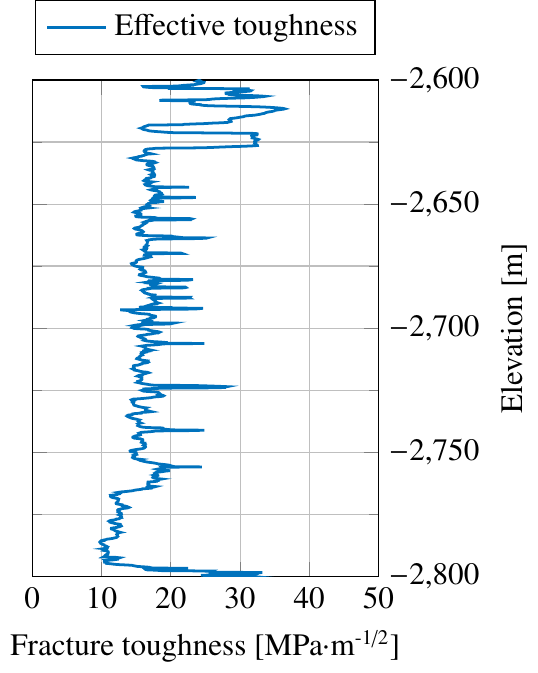}
    \end{subfigure} 
  
  \caption{In-situ rock properties for in the vicinity of the injection.}
  \label{fig:in-situ_profile}
\end{figure*}

We now apply MGR to the simulation of a field-scale 5-stage hydraulic fracturing problem. The computational domain is specified as a 1200$\times$1374$\times$1260~m$^3$ Cartesian block, discretized with 8-node hexahedral elements as shown in \Cref{fig:domain_mesh}.
The computational domain is discretized using permutations of 4 m and 2 m spacing to create 4 refinement levels for the scaling study. These levels correspond to problems of size 21, 42, 84, and 171 million degrees of freedom, respectively.
A horizontal well runs parallel to the y-axis of the problem but is not explicitly modeled here.
The injection points are located at an approximate depth of $-2733$ m, and have an intra-stage spacing of 18 m along the wellbore, with an inter-stage spacing of 36 m.
The in-situ stress profile, elastic properties, and effective fracture toughness estimates are taken from the geologic model of a real unconventional reservoir, and contain extensive vertical heterogeneity as shown in \Cref{fig:in-situ_profile}.
The minimum principal horizontal stress is oriented in the y-direction, while the maximum principal horizontal stress is oriented in the x-direction.

It is important to note that the term effective toughness is an estimate of the field scale toughness including the effects of heterogeneity below the scale of the grid resolution.
The injection fluid is specified with a dynamic viscosity of $0.01$ Pa$\cdot$s (10$\times$ water) and a bulk modulus of $2$ GPa.
The boundary conditions are specified as fixed displacement surrounding the outside of the computational domain, with a mass injection rate of approximately 26.5 kg/s (10 bpm) injected at each perforation location.
Additionally, because the goal of this work is to illustrate the efficacy of the proposed MGR method, we inject at all stages simultaneously in order to activate the coupled flow and mechanics in a larger portion of the computational domain.

\begin{figure}[b!]
    \small
    \centering
    \begin{subfigure}[b]{0.31\textwidth}
      \centering
      \includegraphics[scale=1]{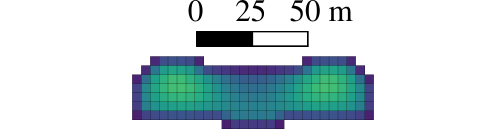}
      \caption{T = 308s.}
    \end{subfigure}
    \hfill
    \begin{subfigure}[b]{0.31\textwidth}
      \centering
      \includegraphics[scale=1]{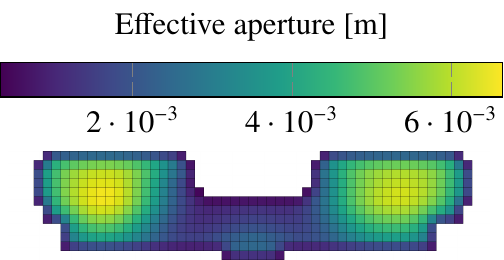}
      \caption{T = 908s.}
    \end{subfigure} 
    \hfill
    \begin{subfigure}[b]{0.31\textwidth}
      \centering
      \includegraphics[scale=1]{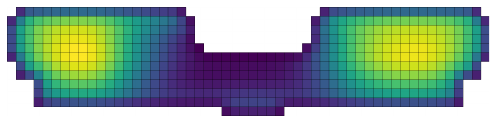}
      \caption{T = 1200s.}
    \end{subfigure} 
    \label{fig:single_fracture_aperture_pressure}
    \caption{Effective aperture of the first fracture in the first stage at different simulation time T(s).}
\end{figure}

\Cref{fig:all_fractures_aperture_full,fig:all_fractures_aperture_zoom} display the aperture of all the fractures embedded in the solid matrix mesh at simulation time T = 908s. 
In \Cref{fig:single_fracture_aperture_pressure}, we show the aperture in the first fracture of the first stage at various times. The effective aperture profiles are well aligned with the in-situ stress shown in \Cref{fig:in-situ_profile}. 
In \Cref{tab:mgr_comparison}, we compare the average number of linear iterations per nonlinear iteration, the setup and solve time of the two preconditioning strategies produced by MGR. 
We use a linear tolerance of $10^{-4}$, which is the same as the nonlinear tolerance. 
Here, although the setup time is comparable between the two methods, MGR\_U is clearly more effective than MGR\_P in terms of convergence, requiring 40\% fewer iterations on average. 
The solve time is also reduced by approximately 40\% as a result. 
One possible reason why MGR\_P performs poorly is because the coarse-grid solve uses an unknown-based approach, which ignores the off-diagonal coupling blocks between $u_x, u_y, u_z$ except for the coarsest level. 
Thus, important off-diagonal values created by reduction with $RAP$ are not taken into account in the coarse-grid solve, which could degrade the convergence.

\begin{table}[t]
    \small
    \centering
    \caption{Comparison between two MGR reduction strategies after the first fracturing stage with linear tolerance \textit{tol} = $10^{-4}$. MGR\_X, X=\{U,P\} denotes choosing the variable X as the F-points for the reduction step.}
    \begin{tabular}{lcc}
        \toprule
        Method & MGR\_U & MGR\_P \\
        \midrule
        Number of cores & 216 & 216 \\
        Total problem size & $\sim$ 21,000,000 & $\sim$ 21,000,000 \\
        Avg. problem size per core & $\sim$ 100,000 & $\sim$ 100,000 \\
        Avg. number of iterations & 12.66 & 21.78 \\
        Avg. setup time per linear solve (s) & 1.33 & 1.31 \\
        Avg. solve time per linear solve (s) & 1.94 & 3.29 \\
        \bottomrule
    \end{tabular}
    \label{tab:mgr_comparison}
\end{table}
\par Based on the results in \Cref{tab:mgr_comparison}, we perform a full simulation run for all the five fracturing stages using MGR\_U. 
The simulation requires approximately 4 times more nonlinear iterations and run time, compared to the first-stage-only simulation. 
The results are shown in \Cref{tab:mgr_u}. It is clear that MGR\_U is robust with respect to the size of the $A_{pp}$ block. 
Even though the fractures expand significantly, and the fluid pressure becomes more prominent, the number of linear iterations using MGR\_U remains stable. 
The total time spent in the linear solve, including both setup and solve phases, still takes up approximately 46\% of the total simulation time, which is quite significant. 
Note that for a less efficient solver, however, this ratio would be much higher. 
Additionally, since data reuse is enabled for MGR\_U, the total time spent in the setup phase is less than 8\% of the total linear solve time.
\begin{table}[t]
    \small
    \centering
    \caption{Performance of MGR\_U for the full hydraulic fracturing simulation with 21M degrees of freedom using 216 cores and \textit{tol} = $10^{-4}$.}

    \begin{tabular}{lc}
        \toprule
        Method & MGR\_U \\
        \midrule
        Avg. number of iterations & 13.00 \\
        Avg. setup time per linear solve (s) & 1.27 \\
        Avg. solve time per linear solve (s) & 2.00 \\
        Total MGR setup time (s) & 967 \\
        Total linear solve time (s) & 12,530 \\
        Total simulation time (s) & 27,455\\
        \bottomrule
    \end{tabular}
    \label{tab:mgr_u}
\end{table}
\par To further challenge the linear solver and examine its weak scalability, we increase the linear tolerance to $10^{-6}$ and perform a series of mesh refinements. 
We successively refine in each direction $x,y,z$ by a factor of 2, resulting in an 8 times increase in the number of degrees of freedom for the largest problem compared to the original one. 
We also increase the number of processor cores for each refinement, keeping the problem size per processor core constant. 
This means that for the largest problem of size 171 million of unknowns, we use 1728 cores on 48 nodes of Quartz. 
The results for MGR\_U are reported in \Cref{tab:weak_scaling}.
\par Since there is strong heterogeneity in the vertical direction of the stress profile, only refining the $x$ and $y$ direction magnifies anisotropy of the problem. 
This makes it more difficult to solve the mechanics block $A_{uu}$ with AMG, which is manifested in a slight growth in number of iterations for the problems with 42M and 84M degrees of freedom. 
When there is uniform refinement in all directions, we achieve a near constant number of linear iterations, despite an 8-fold increase in mesh size. 
Overall, the algorithm shows mesh independent convergence as desired. 
In terms of execution time, the setup time experiences some growth, especially for the largest problem. 
The total setup cost is dominated by the AMG setup of the mechanics block $A_{uu}$. 
Since we are using an aggressive coarsening strategy for the first level and HMIS coarsening for the subsequent levels, both of which uses long-range interpolation operators that involve complicated and costly communication patterns, such an increase is expected and the result is well within known experimental bounds \cite{Baker12}. 
Additionally, the separate displacement component approximation of $A_{uu}$ significantly sparsifies this block and reduces the computational intensity of the setup phase as a result. 
The growth in solve time is much smaller than that of the setup time, increasing only 22\% for the largest problem, compared to the original one. 
This is more important, as the setup time only takes about 1/10th of the solve time if we enable data reuse for MGR\_U. 
Overall, the results are very satisfactory even for large core counts. 
This shows that MGR provides a good platform for rapid discovery of scalable precondtioners.
\begin{table}[ht]
    \small
    \centering
    \caption{Weak scaling for MGR\_U after 320s of hydraulic fracturing simulation with various mesh sizes and \textit{tol} = $10^{-6}$.}
    \begin{tabular}{lcccc}
        \toprule
        Problem size & $\sim$21,000,000 & $\sim$42,000,000 & $\sim$84,000,000 & $\sim$171,000,000\\
        \midrule
        Number of cores & 216 & 432 & 864 & 1,728 \\
        Avg. problem size per core & $\sim$ 100,000 & $\sim$ 100,000 & $\sim$ 100,000 & $\sim$ 100,000\\
        Avg. number of iterations & 29.4 & 35.2 & 38.8 & 30.8 \\
        Avg. setup time per linear solve (s) & 1.23 & 1.72 & 1.70 & 2.65\\
        Avg. solve time per linear solve (s) & 4.35 & 5.51 & 6.35 & 5.33\\
        \bottomrule
    \end{tabular}
    \label{tab:weak_scaling}
\end{table}

\begin{table}[ht]
    \small
    \centering
    \caption{Strong scaling for MGR\_U after 320s for the hydraulic fracturing problem with 21 million degrees of freedom and \textit{tol} = $10^{-6}$.}
    \begin{tabular}{lccccc}
        \toprule
        Number of cores & 108 & 216 & 432 & 864 & 1728 \\
        \midrule
        Avg. problem size per core & $\sim$200,000 & $\sim$100,000 & $\sim$50,000 & $\sim$25,000 & $\sim$12,500\\
        Avg. number of iterations & 27.5 & 29.4 & 28.0 & 29.4 & 28.3\\
        Avg. setup time (s) & 2.20 & 1.23 & 1.28 & 1.28 & 1.45\\
        Avg. solve time (s) & 8.12 & 4.35 & 2.19 & 1.24 & 0.69\\
        \bottomrule
    \end{tabular}
    \label{tab:strong_scaling}
\end{table}

\par Regarding strong scaling, we run the problem with 21M unknowns using 108 cores up to 1728 cores, doubling the number of cores at each stage. 
As can be seen in \Cref{tab:strong_scaling}, MGR scales well. 
The number of iterations remains constant with increasing number of cores. 
In terms of execution time, the solve phase, which largely consists of matrix-vector multiplies, also scales well as expected, almost halving the time for each stage except for a slight loss of efficiency when the average problem size per core becomes small at about 12.5K. 
The setup phase, in contrast, reaches the strong scaling limit much earlier at 100K unknowns per core. 
Below that, the setup time stops decreasing for larger number of cores. 
As explained in the weak scaling section, the setup phase is communication-bound due to the choice of coarsening strategy and approximation of the $A_{uu}$ operator. 
Again, this is a minor issue, as the setup time is less than 10\% of the solve time, and increasing the number of cores still helps achieve faster time to solution.

%% file: Section6_Conclusions.tex
\section{Conclusion and Future Work} \label{sec:conclusions}
In this paper, we present a preconditioning framework based on multigrid reduction to construct efficient and scalable preconditioners for block linear systems that come from simulations of single-phase flow using a hybrid mimetic finite difference discretization, compositional flow with wells, and hydraulic fracturing. Mesh independent convergence and strong scalability results are obtained for the compositional flow with wells and single-phase flow using a hybrid mimetic finite difference method for field-scale problems with a heterogeneous permeability field. In the hydraulic fracturing case, the MGR preconditioner is able to capture the strong coupling between fracture mechanics and pressure-driven fluid flow, resulting in an efficient method that exhibits mesh independent convergence and scales well to thousands of cores and problems with hundred millions of unknowns. To our knowledge, this is the first time that an algebraic solver can achieve performance on this scale for this class of problem. 
\par For future work, we would like to combine hydraulic fracturing with compositional flow by considering compositional displacements inside the fracture and adding fracture-matrix fluid exchange. We will explore compositional flow with a large number of components, and will generalize the hybrid MFD linear solution strategy to compositional flow with wells. These extensions will result in larger and more complex block systems with very different types of couplings. The MGR framework developed here would provide exactly the tool for rapid development of efficient preconditioners that can tackle such complex systems.


%% file: Appendix.tex
\crefalias{section}{appendix}

\section{Hybrid Mimetic Finite Difference for Single-phase Flow}\label{app:hybrid_mfvm}

In Section~\ref{sec:results}, compressible single-phase flow is governed by the following system consisting of Darcy's law and a mass balance equation:
\begin{linenomath}
\begin{subequations}
\begin{align}
& \boldsymbol{u} = - \frac{\kappa}{\mu} ( \nabla p - \rho \boldsymbol{g} ), & & \mbox{(Darcy's law),} \label{eq:single_phase_darcy} \\
& \frac{\partial}{\partial t} \left( \phi \rho \right) + \nabla \cdot \left( \rho \boldsymbol{u} \right) = 0, & & \mbox{(mass balance).} \label{eq:single_phase_mass_balance}
\end{align} \label{eq:IBVP_single_phase}\null
\end{subequations}
\end{linenomath}
where $\boldsymbol{u}$ is the Darcy velocity, $\kappa$ is the (scalar) permeability coefficient, $p$ is the pressure, and $\boldsymbol{g}$ is the gravity vector weighted by gravitational acceleration.
The porosity, $\phi$, the viscosity, $\mu$, and the density, $\rho$, depend linearly on pressure.
This system is discretized with a hybrid MFD scheme described in \cite[Chap. 6]{Lie19}.
This formulation involves three types of degrees of freedom (dofs).
In this appendix only, we use a bold straight font to distinguish the vectors of dofs from the unknowns of \Cref{eq:IBVP_single_phase}.
The set of cell-centered dofs, $\mathbf{p}$, contains cell-based pressure averages.
The set of one-sided face dofs, $\mathbf{w}$, is a collection of face-based flux averages.
In the hybrid formulation, there is one one-sided flux dof per boundary face, and two one-sided flux dofs per interior face (not necessarily equal).
The set of Lagrange multipliers, $\boldsymbol{\pi}$, can be interpreted as a collection of face-based pressure averages.
There is one multiplier per face.

The mass balance equation and Darcy's law are assembled as separate (coupled) equations into the Jacobian matrix.
\Red{In Sections \ref{subsec:hybrid_mfvm} and \ref{subsec:hybrid_mfvm_results}, we consider two discretizations of \Cref{eq:single_phase_darcy}, referred to as the TPFA approach and the consistent approach.
In the TPFA hybrid MFD approach, \Cref{eq:single_phase_darcy} is discretized using the standard Two-Point Flux Approximation scheme.
In the consistent hybrid MFD approach, a consistent discretization of \Cref{eq:single_phase_darcy} on non K-orthogonal meshes is obtained with the introduction of a coercive discrete bilinear form referred to as an inner product.
We use here the (consistent) parametric inner product defined in \cite[Chap. 6]{Lie19} and set the parameter multiplying the stability term to 2.
Both discretizations can be written in the form:}
\begin{linenomath}
\begin{equation}
A_{ww} \mathbf{w} - D \mathbf{p} - C \boldsymbol{\pi} - \mathbf{b}(\mathbf{p}) = \mathbf{0}. \label{eq:discrete_darcy}
\end{equation}
\end{linenomath}
Here, $A_{ww}$ is a block-diagonal SPD matrix containing one block per cell.
\Red{Each block is diagonal in the TPFA approach and dense in the consistent approach.}
$D$ is a rectangular matrix whose entry $[D]_{ij}$ is equal to one if one-sided face $i$ belongs to cell $j$, and zero otherwise.
The linear operator $\mathbf{b}(\mathbf{p})$ takes the vector of cell-centered pressure dofs, $\mathbf{p}$, in input, and outputs a vector of gravity terms in the space of one-sided face fluxes. %
Specifically, we have $[\mathbf{b}(\mathbf{p})]_{i} = - \rho_j g_{ij}$ where $g_{ij}$ is obtained by multiplying the gravitational acceleration by the depth difference by the center of cell $j$ and the center of one-sided face $i$.
%
%
$C$ is a rectangular matrix whose entry $[C]_{ij}$ is equal to -1 if one-sided face $i$ belongs to face $j$, and zero otherwise.

The discrete mass balance equation involves the accumulation term and the divergence of the one-sided fluxes weighted by density. It reads:
\begin{linenomath}
\begin{equation}
\mathbf{m} - \mathbf{m}^n + \Delta t D^T( U( \mathbf{p} ) \mathbf{w} ) = \mathbf{0}, \label{eq:discrete_mass_balance}
\end{equation}
\end{linenomath}
where $\mathbf{m}$ denotes a vector of cell masses written as the product of volume, porosity, and density such that $[\mathbf{m}]_i = v_i \phi_i \rho_i$. 
In the diagonal matrix $U( \mathbf{p} )$, the $i$-th diagonal entry contains the upwinded density at one-sided face $i$.
Note that wells are not modeled in this example and are therefore absent from this equation.

In the hybrid framework, continuity of the volumetric one-sided fluxes at the mesh faces is enforced through the constraint:
\begin{linenomath}
\begin{equation}
C^T \mathbf{w} = \mathbf{0}. \label{eq:discrete_face_constraints}
\end{equation}
\end{linenomath}
The nonlinear system consisting of \Cref{eq:discrete_darcy,eq:discrete_mass_balance,eq:discrete_face_constraints} is solved with Newton's method.
At each linearization step, setting $A_{wp} = -D - \frac{\partial \mathbf{b}}{\partial \mathbf{p}}(\mathbf{p})$, $A_{w\pi} = C$, $A_{pw} = \Delta t D^T U(\mathbf{p})$, $A_{pp} = \frac{\partial}{\partial \mathbf{p}}( \mathbf{m} + \Delta t D^T( U( \mathbf{p} ) \mathbf{w} ) )$, and $A_{\pi w} = C^T$ yields the block system of \Cref{eq:mfd_saddle_point_system}.

\section{Compositional Flow with Wells}\label{app:compositional_flow}

This section provides an overview of the isothermal two-phase, two-component compositional model used in Section \ref{sec:results} to simulate CO$_2$ injection.
More details about the formulation can be found in \cite{Cao02,Voskov12} and references therein.
The system is governed by a set of mass balance equations---one for each component, namely, H$_2$O and CO$_2$---in which the phase velocities are given by Darcy's law:
\begin{linenomath}
\begin{subequations}
\begin{align}
& \boldsymbol{u}_\ell = - \lambda_\ell \kappa ( \nabla p - \rho_\ell \boldsymbol{g} ), & & \ell = 1, 2, & & \mbox{(Darcy's law),} \label{eq:darcy} \\
& \frac{\partial}{\partial t} \left( \phi \sum_{\ell=1}^{2} x_{c\ell} \rho_\ell s_{\ell} \right) + \nabla \cdot \left( \sum_{\ell=1}^{2} x_{c\ell} \rho_\ell \boldsymbol{u}_\ell \right) + \sum_{\ell=1}^{2} x_{c\ell} \rho_\ell q_\ell = 0, & & c = 1, 2, & & \mbox{(mass balance).} \label{eq:compositional_mass_balance}
\end{align} \label{eq:IBVP_compositional}\null
\end{subequations}
\end{linenomath}
In \Cref{eq:darcy}, $\lambda_{\ell}$ is the phase mobility, $\kappa$ is the (scalar) absolute permeability, $p$ is the pressure, and $\boldsymbol{g}$ is the gravity vector.
In \Cref{eq:compositional_mass_balance}, $\phi$ is the porosity, $x_{c\ell}$ is the phase component fraction, $\rho_\ell$ is the phase (mass) density, $s_\ell$ is the phase saturation representing the fraction of the porous medium occupied by phase $\ell$, and $q_\ell$ is the well term.
We neglect capillary pressure.
The system is closed by the following local constraints:
\begin{linenomath}
\begin{subequations}
\begin{align}
f_{cj} - f_{ck} = 0, \qquad &  \forall j \neq k, \, c = 1, 2, & & \text{(thermodynamic equilibrium),} \label{eq:thermodynamic_equilibrium} \\
x_{1\ell} + x_{2\ell} - 1 = 0, \qquad &  \ell = 1, 2, & & \text{(component fraction constraints),} \label{eq:component_fraction_constraints} \\
s_1+s_2 - 1 = 0, \qquad &  & & \text{(saturation constraint).} \label{eq:saturation_constraint} 
\end{align}
\end{subequations}
\end{linenomath}
The full system consisting of \Cref{eq:compositional_mass_balance,eq:thermodynamic_equilibrium,eq:component_fraction_constraints,eq:saturation_constraint} involves $n_p(n_c+1)+1 = 7$ equations and unknowns.
The standard solution strategy consists in identifying a set of primary equations that are differentiated with respect to the primary unknowns and assembled into the Jacobian matrix, while secondary equations and unknowns are eliminated locally during the assembly.

The three primary equations considered here are the two mass balance equations (\Cref{eq:compositional_mass_balance}) and the saturation constraint (\Cref{eq:saturation_constraint}).
The three primary unknowns are the pressure, $p$, and the two overall component densities, $\rho_c$ ($c = 1, 2$).
An overall component density represents the mass of a given component per unit volume of mixture, and can be related to the variables of \Cref{eq:compositional_mass_balance} using the formulas given in \cite{Voskov12}.
The partitioning of the CO$_2$ component between the CO$_2$ phase and the aqueous phase is determined using the model of Duan and Sun \cite{Duan03}.
The CO$_2$ phase densities and viscosities are computed using the Span-Wagner \cite{Span96} and Fenghour \cite{Fenghour98} correlations, respectively, while the brine properties are obtained using the methodology of Phillips et al. \cite{Phillips81}.

The system is discretized with a cell-centered, fully implicit, finite-volume scheme based on single-point upstream weighting \cite{Cao02}, and solved with Newton's method with damping.
At each linearization step, we obtain a linear system with the structure given in \Cref{eq:compositional_reservoir_block}, in which the mass balance equation for CO$_2$ is aligned with the pressure, the mass balance equation for H$_2$O is aligned with component density of CO$_2$, and the saturation constraint is aligned with the component density of H$_2$O.
The well flow formulation relies on the same degrees of freedom, to which we add a degree of freedom representing the total mass rate flowing between well segments.
The well transmissibility at the perforations is computed using the Peaceman formula.

\section{Hydraulic Fracturing} \label{app:hydraulic_fracturing}
In this section, we provide an overview of the governing equations and discretization strategy used within our hydraulic fracturing simulator.  A more extensive description can be found in \cite{Settgast:2017}.  As the main focus of the work is the linear solver strategy, we limit the presentation here to the essential details.  For testing purposes, we also focus on a simplified model that excludes some physics---notably fluid leakoff---that is important in practical applications but that can be readily addressed once a baseline solver strategy is in place.

We consider a deformable elastic body $\mathcal{B}$ that is cut by one or more internal fracture surfaces $\mathcal{F}$.  Fluid is injected into these fractures through a well.  Due to the fluid pressure, the elastic body will deform so that the aperture of the fractures grows and the space available for fluid storage increases.  This deformation also causes a strong stress concentration at the outer edges of the fracture.  Once the stress intensity at an edge exceeds a certain threshold---a mechanical property of the rock---the rock will break and the fracture can propagate.  We assume for simplicity that the elastic body is impermeable, so that all fluid remains within the fracture(s).  

The simulation of this process is broken into two key steps, which are used iteratively to evolve the system in time.  Let $t_n$ denote the current time at timestep $n$, and let $k$ denote an iteration counter.  Beginning from an initial guess at the configuration of the system, $k=0$, the following steps are repeated until convergence:
\begin{enumerate}
\item \textit{Equilibration}: We assume the boundary of the fractures $\partial \mathcal{F}^k$ is known.  We solve governing mass and momentum balance equations to determine a new displacement of the elastic body $\vec{u}^k_n$ and fluid pressure $p^k_n$ in the fracture.
\item \textit{Propagation}: Using the updated state, stress conditions at the edge of the fracture can be evaluated and a new fracture boundary $\partial \mathcal{F}^{k+1}$ determined.  Convergence is reached when the boundary stabilizes, $\partial \mathcal{F}^{k} = \partial \mathcal{F}^{k+1}$.  Otherwise, we set $k \gets k+1$ and a new equilibration solve is required.
\end{enumerate}
The domain $\mathcal{B}$ is discretized using a triangulation that internally conforms to the fracture surface $\mathcal{F}$. As the fracture propagates, faces within this triangulation are ``broken'' and the mesh topology recomputed.  This propagation algorithm is much cheaper than the equilibration solve and highly localized.   

In contrast, during equilibration a momentum balance equation is solved on the volumetric mesh to determine the deformation of the elastic body, and a mass balance equation is solved on a lower dimensional mesh made up of the union of broken fracture faces.  The two equations are strongly coupled to each other through several interaction terms.  In particular, we seek $(\vec{u},p)$ such that
\begin{linenomath}
\begin{align}
\nabla \cdot \vec{\sigma} + \rho \vec{g} = \vec{0} &\quad \text{on} \; \mathcal{B} \setminus \mathcal{F} \quad &  & \text{(momentum balance),} \\
\dot{m} + \nabla \cdot \vec{q} = s &\quad \text{on} \; \mathcal{F} \quad & &\text{(mass balance),} 
\intertext{subject to the internal boundary condition}
\vec{\sigma} \cdot \vec{n} = p \vec{n}  &\quad \text{on} \; \mathcal{F} \quad  & &\text{(traction balance).}
\end{align}
\end{linenomath}
Here, the stress $\vec{\sigma} = \mathbb{C}:\nabla^s \vec{u}$, where $\mathbb{C}$ is a fourth-order elasticity tensor and $\nabla^s$ is the symmetric gradient; $\rho$ is the solid density and $\vec{g}$ is gravity; $\vec{n}$ is the local unit normal the fracture surface; $m = \rho_f w$ is the fluid mass per unit fracture area, which is determined from the fluid density $\rho_f$ and fracture aperture $w = \llbracket \vec{u} \rrbracket \cdot \vec{n}$, where $\llbracket \vec{u} \rrbracket$ indicates the jump in the displacement field across the fracture; $\vec{q} = - \rho_f (w^3 / 12 \eta) ( \nabla p - \rho_f \vec{g})$ is the fluid mass flux through the fracture, with fluid viscosity $\eta$; and $s$ is a fluid source term representing injection wells.  We remark that the mass flux through the fracture is a strongly nonlinear function of the aperture.  These governing equations are supplemented with appropriate initial conditions and external boundary conditions.